\documentclass[letterpaper, 10 pt, conference]{ieeeconf}
\usepackage{ifpdf, flushend,subfigure}
\usepackage{graphicx}
\usepackage{epstopdf}
\graphicspath{{../}}
\usepackage{multirow}
\usepackage{float}
\usepackage[cmex10]{amsmath}
\usepackage {amssymb}
\usepackage{array}
\usepackage{mdwmath}
\usepackage{mdwtab}
\usepackage{eqparbox}
\usepackage{url}
\usepackage{nomencl}
\usepackage{cite}
\usepackage{algorithm}
\usepackage{algorithmic}
\usepackage{makeidx}
\usepackage{ifthen}
\usepackage{subfigure}
\usepackage{caption}
\usepackage{pdflscape}
\usepackage{hyperref}
\usepackage{xcolor}
\def\mbi{\mathbb{I}}

\def\rev#1{{\color{black}#1}} 

\newcommand{\FE}{{\textnormal{FE}}}
\newcommand{\FR}{{\textnormal{FR}}}
\newcommand{\FB}{{\textnormal{FB}}}
\newcommand{\PG}{{\textnormal{PG}}}
\newcommand{\PR}{{\textnormal{PR}}}
\newcommand{\PS}{{\textnormal{PS}}}
\newcommand{\PC}{{\textnormal{PC}}}
\newcommand{\PB}{{\textnormal{PB}}}
\newcommand{\PD}{{\textnormal{PD}}}

\makeatletter
\def\munderbar#1{\underline{\sbox\tw@{$#1$}\dp\tw@\z@\box\tw@}}
\makeatother

\usepackage{amsmath,mathtools, amssymb, bbm, xspace}
\usepackage{booktabs,ragged2e}
\usepackage{multirow}
\usepackage[flushleft]{threeparttable}

\usepackage{graphicx}

\def\det{\mathop{\hbox{\rm det}}}

\def\spose#1{\hbox to 0pt{#1\hss}}

\def\text #1{\hbox{\quad#1\quad}}

\def\nthinsp{\mskip -2   mu}

\def\superstar{^{\raise 0.5pt\hbox{$\nthinsp *$}}}
\def\SUPERSTAR{^{\raise 0.5pt\hbox{$*$}}}

\def\lamstarT {\lambda^{\raise 0.5pt\hbox{$\nthinsp *$}T}}

\def\hbar{\skew{4.2}\bar h}

		\def\bkE{{\rm I\kern-.17em E}}
		\def\bk1{{\rm 1\kern-.17em l}}
		\def\bkD{{\rm I\kern-.17em D}}
		\def\bkR{{\rm I\kern-.17em R}}
		\def\bkP{{\rm I\kern-.17em P}}
		\def\bkY{{\bf \kern-.17em Y}}
		\def\bkZ{{\bf \kern-.17em Z}}

		\def\bc{\begin{center}}
		\def\be{\begin{enumerate}}
		\def\bi{\begin{itemize}}
		\def\ec{\end{center}}
		\def\ee{\end{enumerate}}
		\def\ei{\end{itemize}}
		\def\es{\end{small}}
		\def\eS{\end{slide}}

	\def\cp2problem#1#2#3#4{\fbox
		 {\begin{tabular*}{0.9\textwidth}
			{@{}l@{\extracolsep{\fill}}l@{\extracolsep{6pt}}l@{\extracolsep{\fill}}c@{}}
				#1 & & $#4 $ 
			\end{tabular*}}}

		\def\bkE{{\rm I\kern-.17em E}}
		\def\bk1{{\rm 1\kern-.17em l}}
		\def\bkD{{\rm I\kern-.17em D}}
		\def\bkR{{\rm I\kern-.17em R}}
		\def\bkP{{\rm I\kern-.17em P}}
		
		\def\bkZ{{\bf{Z}}}

\newcommand {\beeq}[1]{\begin{equation}\label{#1}}
\newcommand {\eeeq}{\end{equation}}
\newcommand {\bea}{\begin{eqnarray}}
\newcommand {\eea}{\end{eqnarray}}

\def\texitem#1{\par\smallskip\noindent\hangindent 25pt
               \hbox to 25pt {\hss #1 ~}\ignorespaces}

\newcommand{\beq}{\begin{equation}}
\newcommand{\eeq}{\end{equation}}
\newcommand{\beqn}{\begin{eqnarray}}
\newcommand{\eeqn}{\end{eqnarray}}
\newcommand{\beqno}{\begin{eqnarray*}}
\newcommand{\eeqno}{\end{eqnarray*}}
\newcommand{\bma}{\begin{displaymath}}
\newcommand{\ema}{\end{displaymath}}
\newcommand{\bnu}{\begin{enumerate}}
\newcommand{\enu}{\end{enumerate}}
\newcommand{\bce}{\begin{center}}
\newcommand{\ece}{\end{center}}
\newcommand{\btb}{\begin{tabular}}
\newcommand{\etb}{\end{tabular}}

\newcommand{\col}{\mathrm{col}}
\def\R{{\mathbb{R}}}
\def\G{{\mathcal{G}}}

\def\bD{{\mathbf{D}}}
\def\bS{{\mathbf{S}}}
\def\R{{\mathbf{R}}}
\def\C{{\mathbf{C}}}
\def\E{{\mathcal{E}}}

\def\bx{{\boldsymbol{x}}}

\def\by{{\boldsymbol{y}}}
\def\bz{{\boldsymbol{z}}}
\def\bu{{\boldsymbol{u}}}
\def\bv{{\boldsymbol{v}}}

\def\bQ{{\mathbf{Q}}}
\def\bG{{\mathbf{G}}}
\def\bY{{\mathbf{Y}}}

\def\b1{{\mathbf{1}}}

\usepackage[normalem]{ulem} 

\newcommand{\mc}{\mathcal}

\newtheorem{theorem}{Theorem}
\newtheorem{definition}{Definition}

\newtheorem{lemma}{Lemma}
\newtheorem{remark}{Remark}

\newtheorem{assumption}{Assumption}

\newcommand{\g}{{\mathbf{g}}}

\newcommand{\Diag}{{\mathrm{Diag}}}

\newcommand{\ra}{{r_\mathrm{a}}}

\newcommand{\V}{{\mathcal{V}}}

\newcommand{\one}{{\mathbf{1}}}
\newcommand{\zero}{{\mathbf{0}}}

\newcommand{\bE}{{\mathbf{E}}}

\newcommand{\bw}{{\boldsymbol{w}}}

\newcommand{\T}{{\mathsf{T}}}

\definecolor{myBlue}{rgb}{0.80,0.85,1.00}
\definecolor{myYellow}{rgb}{0.951,1.000,0.547}

\def\la{{\langle}}
\def\ra{{\rangle}}

\def\bit{\begin{itemize}}
\def\eit{\end{itemize}}
\def\BEAS{\begin{eqnarray*}}
\def\EEAS{\end{eqnarray*}}

\def\re{{\mathbb R}}

\def\a{\alpha}
\def\b{\beta}

\def\g{\gamma}

\IEEEoverridecommandlockouts     
\title{\LARGE \bf
Constrained multi-cluster game: Distributed Nash equilibrium seeking over directed graphs
}

\author{Duong Thuy Anh Nguyen\textsuperscript{1}, Mattia Bianchi\textsuperscript{2}, Florian Dörfler\textsuperscript{2},
Duong Tung Nguyen\textsuperscript{1}, Angelia Nedi\'c\textsuperscript{1}
\thanks{1. School of Electrical, Computer and Energy Engineering, Arizona State University, Tempe, AZ, United States. 
Email: \{dtnguy52,~duongnt,~Angelia.Nedich\}@asu.edu. This work was supported in parts by the NSF award CCF-2106336, the ONR award N00014-21-1-2242, and the ARPA-H award SP4701-23-C-0074. 
2. Automatic control laboratory, ETH Zürich, Switzerland. Email: \{mbianch,dorfler\}@ethz.ch. This work was supported by ETH Zürich funds. 
}
}

\begin{document}
\maketitle

\begin{abstract} 
Motivated by the complex dynamics of cooperative and competitive interactions within networked agent systems, multi-cluster games provide a framework for modeling the interconnected goals of self-interested clusters of agents. For this setup, the existing literature lacks comprehensive gradient-based solutions that simultaneously consider constraint sets and directed communication networks, both of which are crucial for many practical applications. To address this gap, this paper proposes a distributed Nash equilibrium seeking algorithm that integrates consensus-based methods and gradient-tracking techniques, where inter-cluster and intra-cluster communications only use row- and column-stochastic weight matrices, respectively. To handle constraints, we introduce an averaging procedure, which can effectively address the complications associated with projections. In turn, we can show linear convergence of our algorithm, focusing on the contraction property of the optimality gap. We demonstrate the efficacy of the proposed algorithm through a microgrid energy management application.

\end{abstract}

\printnomenclature
\allowdisplaybreaks

\section{Introduction} \label{intro}

In various networked systems across different domains, such as telecommunications \cite{bianchi2024end}, transportation systems \cite{Acciaio2021}, crowdsourcing \cite{wiopt23} and energy distribution \cite{Liu2021,BIANCHI2021109660}, 
interconnected entities operate autonomously, with both cooperative and competitive emerging strategies.
To capture these diverse dynamics,
multi-cluster games \cite{YE2018266} emerge as a pertinent framework, expanding traditional non-cooperative game \cite{BIANCHI2022110080} and distributed optimization setting \cite{Nedic2009} to scenarios where clusters of agents compete. 
In multi-cluster games, each cluster represent a cohesive group of agents collaborating to minimize the sum of their local objective functions. Yet, clusters compete with each other  
in a noncooperative game framework, independently making decisions to achieve distinct but interconnected goals. For instance, in smart grids, multiple microgrids compete in energy distribution, while at the local level, an economic dispatch problem is solved. 
Nash equilibria (NEs) serve as a pivotal tool for analyzing and optimizing decision-making in such decentralized systems.

In recent years, significant interest has emerged in algorithms for  NE seeking \cite{Bianchi_Directed_CDC_2020,Bianchi_TV_LCSS_2021,Tatarenko2021,Nguyen2023RowGame,nguyen2022distributed,Hall2022} and for distributed optimization \cite{Tsianos2012,nedic2015distributed,nedic2017achieving,xi2018add,Qu_Harnessing_TCNS2018,pshi21,Angelia2022AB,Nguyen2023AccAB,Nguyen2023SABTV}, particularly under the partial-decision information scenario. In such settings, agents are constrained to access solely their own cost functions and local action sets, with limited information exchange with neighboring agents on a communication network. Motivated by these advancements, numerous distributed algorithms
\cite{YE2018266,Nian2022,YE2020108815,ZENG201920,Sun2021,MENG2023110919,Zimmermann2022,Zimmermann2021,pang2023distributed,Zhou2023Multicoalition} have  also  been proposed for multi-cluster games, leveraging local information exchange among agents.

A substantial body of research has been dedicated to continuous-time algorithms for finding NEs\cite{YE2018266,Nian2022,YE2020108815}, and generalized Nash equilibria (GNEs) in multi-cluster games \cite{ZENG201920,Sun2021}.
The work in \cite{YE2018266} solves  unconstrained games using gradient-based algorithms, while \cite{Nian2022} extends this to directed time-varying communication topologies. 
Concerning the presence of constraints, reference \cite{ZENG201920} investigates a distributed  projected differential inclusion, to find GNEs in nonsmooth games with coupled nonlinear inequality constraints and set constraints. Instead, \cite{Sun2021} addresses games with inequality constraints, employing finite-time average consensus. Both papers only consider undirected graphs.

For discrete-time algorithms,
some works focus on gradient-free and payoff-based methods \cite{Tatarenko2021ClusterGames}. In contrast to these solutions,   
gradient-based algorithms do not assume that the agents can measure their costs, but require inter-cluster communication, as it is necessary for the agents to estimate the joint strategy in order to evaluate their local gradients. Many proposed algorithms are inspired by consensus-based algorithms originally developed for distributed optimization \cite{Qu_Harnessing_TCNS2018,pshi21,Angelia2022AB,Nguyen2023AccAB,Nguyen2023SABTV}, incorporating gradient-tracking techniques by limiting gradient tracking to agents within the same cluster. In \cite{MENG2023110919}, a leader-follower hierarchy is established, where followers within each cluster communicate solely with their in-cluster neighbors and leader, and inter-cluster communication is limited to exchanges between cluster leaders. The work in \cite{Zimmermann2022} extends these findings by introducing a more general leaderless communication architecture.
All communications in these studies are undirected. References \cite{Zimmermann2021,Zhou2023Multicoalition,pang2023distributed} propose similar algorithms, but with a directed communication network -- yet, their results are limited to unconstrained games.

While there exist several versions of gradient-based algorithms for distributed NE seeking in multi-cluster games \cite{YE2018266,YE2020108815,ZENG201920,Sun2021,MENG2023110919,Zimmermann2022,pang2023distributed,Zimmermann2021,Zhou2023Multicoalition}, none of them addresses the presence of constraints and directed communication simultaneously. Incorporating both is crucial for many practical applications \cite{Liu2021,Chengcheng2017,Zimmermann2021,Zimmermann2022}, 
where constraints represent real-world limitations such as physical capacities, operational requirements, and regulatory standards. Yet, employing projection methods is not straightforward: indeed, the applicability of gradient-tracking techniques alongside a proximal method \cite{FALSONE2022109938} in multi-cluster games remains an open issue. This challenge is 
further exacerbated by directed communications, which introduce technical complexities stemming from
imbalanced weight matrices.
{Yet, dealing with directed graphs is essential, 
as unilateral communication capabilities arise naturally in wireless scenarios, for instance when communication ranges of different agents (e.g., sensors, smart meters) vary due to heterogeneous wireless transmission technologies, differing transmitting power levels, or fluctuating channel conditions/noises. Communication interference  and cyberattacks can also result in directed communication networks.}

\textbf{Contributions.} Drawing from the works \cite{Nguyen2023RowGame} for non-cooperative games and \cite{Orhan2023} for distributed optimization, we propose a distributed NE seeking algorithm for multi-cluster games that incorporates the consensus-based approach, the gradient-tracking technique, and an averaging procedure (or ``lazy update''), effectively addressing  technical complications associated with constraint sets and directed communication networks. Particularly, choosing  the averaging parameter appropriately allows us to control the error in the updates generated by applying projection methods within the gradient-tracking procedure, and to ensure linear convergence of our algorithm. We numerically evaluate our method on a microgrid energy management problem, for which the raw experimental data and the code repository is publicly available\footnote{\textnormal{https://github.com/duongnguyen1601/Distributed\_NCluster\_Game/}}.

We outline the convergence analysis of the proposed algorithm by establishing a contraction relationship among three error terms, for which a detailed analysis of the contraction relation of the optimality gap is provided.
For the contraction relations of the consensus error and the gradient tracking error, we refer to the existing proofs in \cite{Zimmermann2021} and \cite{Orhan2023} 
for unconstrained games and distributed optimization, respectively. The constrained multi-cluster game under consideration and the directed communication pose some challenges when applying the referenced results in our analysis. Therefore, we formulate the updates of the algorithm and the error terms, as well as present relevant lemmas to establish connections with existing results for seamless adaptation to our scenario. 


\textbf{Notations.} 
We let $[n]=\{1,\ldots,n\}$ for an integer $n\ge 1$.
All vectors are column vectors unless otherwise stated. For a vector $u\in \re^n$, we use $u^\T$ to denote its transpose. We define $\min(u)=\min_i u_i$ and $\max(u)=\max_i u_i$. $\Diag(u)$ denotes the diagonal matrix whose diagonal entries correspond to the entries of $u$. A nonnegative vector is called stochastic if its entries sum up to $1$. We use $\zero$ and $\one$ to denote the vector with all entries equal to $0$ and $1$, respectively.

We use $A_{ij}$ to denote the $ij$-th entry of a matrix $A$.
The notation $A\le B$ implies $A_{ij}\le B_{ij}$ for all $i, j$.
A matrix $A$ is nonnegative if all its entries are nonnegative and ${\min}^{+}(A)$ denotes the smallest positive entry of $A$. 
The identity matrix is denoted by $\mbi$. 

Given a vector $\pi\in\re^m$ with positive entries, we denote
\begin{center}
    $\la \bu,\bv\ra_{\pi}=\sum_{i=1}^m \pi_i\la u_i,v_i \ra $ and $\|\bu\|_{\pi}=\sqrt{\sum_{i=1}^m \pi_i\|u_i\|^2},$
\end{center}
where \rev{$\bu\!:=\![u_1,\ldots,u_m]^\T, \bv \!:=\![v_1,\ldots,v_m]^\T \!\!\in\! \re^{m\times n}$}, and $u_i,v_i\!\in\!\re^n$. When $\pi = \one$, we write $\la \bu,\bv\ra$ and $\|\bu\|$. 
We also write $\|\bu\|_{\pi^{-1}}$ to denote the norm induced by the vector with entries $1/\pi_i$, i.e., 
$\|\bu\|_{\pi^{-1}}=\sqrt{\sum_{i=1}^m\|u_i\|^2 /\pi_i}$.

 
A directed graph $\G$  is considered \textit{strongly connected} if there exists a directed path from any node to all other nodes in $\G$. 
We denote the diameter and the maximal edge-utility of a strongly connected directed graph $\G$ as $\mathsf{D}(\G)$ and $\mathsf{K}(\G)$, respectively, as defined in \cite[Definitions 2.1 and 2.2]{Angelia2022AB}.

\section{Problem Formulation} \label{sec:formu}

\subsection{Multi-cluster Game}
Consider the multi-cluster game played by the set $\V$ of $N$ agents, grouped into $H$ clusters. Each cluster operates as a virtual agent within a non-cooperative game framework. Cluster $h\in [H]$ is defined by a subset $\V_h$ of $N_h$ agents, where $N=\sum_{h\in [H]} N_h$. The agent sets for all clusters are disjoint, i.e., $\V_h \bigcap \V_{h'} = \emptyset$ for $h \neq h'$, and $\cup_{h=1}^H \V_h = \V$. 
To simplify the notation, we use superscripts to represent agent indices 
and subscripts to denote cluster indices.

Each agent $i$ within cluster $h$ is associated with a cost function $f^i(x_h,x_{-h})$, known only to agent $i$. This function depends on the action of its own cluster, $x_h\in X_h \subseteq \re^{p_h}$, and the joint action of all other clusters except its own, $x_{-h} \!=\! \col( (x_\ell)_{\ell\in \mc [H]\backslash \{ h \} }\!) \!\in\! X_{-h} \!\subseteq\! \re^{p-p_h}$. The joint action vector of all clusters is $x=\col( (x_h)_{h \in [H]})$, has size $p = \sum_{h\in [H]} p_h$ and belongs to the joint action set $X = X_1\times\cdots\times X_H \subseteq \re^p$. The agents within each cluster $h$ collaborate to minimize the cost function $F_h$, for all $h \in [H]$, as follows:
\begin{align}\label{eq:clusterCostFunction}
\min_{x_h \in X_{h}} F_h(x_h,x_{-h}) = \min_{x_h \in X_{h}} \frac{1}{N_h}\sum_{i\in \V_h} f^i(x_h,x_{-h}).
\end{align}
It is imperative to highlight that agents must achieve consensus on the strategy $x_h$ to minimize the cost function $F_h$ within their respective cluster $h$. In the decision-making process, agents can only adjust the strategy of their own cluster, while observing the strategy of other clusters utilizing information exchange through communication networks.

\begin{remark} 
When the number of clusters is $H=1$, the multi-cluster game is reduced to the distributed optimization problem, as described in \cite{Angelia2022AB}. When the number of agents in each cluster is $N_h = 1$ for all $h\in [H]$, the multi-cluster game reduces to the non-cooperative game of $N$ agents \cite{nguyen2022distributed}.
\end{remark}

Denote the game by $\Gamma=([H],\{F_h\},\{X_h\})$. An NE for the game $\Gamma$ can be formally defined as follows:
\begin{definition}[Nash equilibrium] \label{def-NE}
For the multi-cluster game $\Gamma$, a strategy profile $x^*= \col( (x_h^*)_{h \in [H]})\in X$ is an NE of the game if, for every cluster $h\in [H]$, there holds:
\begin{align*}
F_h(x_h^*,x_{-h}^*)\le F_h(x_h,x_{-h}^*),  \quad\hbox{for all } x_h\in X_h.
\end{align*}
\end{definition}

We define the game mapping $M(\cdot):\re^p\to\re^p$ as follows
\begin{align}\label{eq:gamemapping}
M(x)\triangleq \col\left( (\nabla_h F_h(x_h,x_{-h}))_{h \in [H]}\right),
 \end{align}
where $\nabla_h F_h(x_h,x_{-h})\triangleq \frac{1}{N_h}\sum_{i\in \V_h} \nabla_h f^i(x_h,x_{-h}) \in \re^{p_h}$ with $\nabla_h \triangleq \nabla_{x_h}$, for all $h\in[H]$.

We make the following standard assumptions \cite{Zimmermann2021,MENG2023110919}.

\begin{assumption}\label{assum:lip}
Consider the game $\Gamma$, and assume for all cluster $h \in [H]$ and for all agent $i \in \V_h$:\\
(i) The action set $X_h$ is non-empty, closed and convex.\\ 
(ii) The cost function $f^i(x_h,x_{-h})$ is convex and continuously differentiable in $x_h$ for any fixed $x_{-h}\in \re^{p-p_h}$.\\
(iii) The gradient $\nabla_h f^i(x_h,x_{-h})$ is Lipschitz continuous on
$\re^{p-p_h}$ for every fixed $x_h\in \re^{p_h}$ with a constant $L_1^i>0$.\\
(iv) The gradient $\nabla_h f^i(x_h,x_{-h})$ is Lipschitz continuous on $\re^{p_h}$ for every fixed $x_{-h}\in \re^{p-p_h}$ with a constant $L_2^i>0$.
\end{assumption}

Given Assumption~\ref{assum:lip}(i)-(ii), an NE of the game $\Gamma$ can alternatively be characterized through the first-order optimality conditions. Specifically, $x^*\in X$ is an NE of the game $\Gamma$ if and only if,
for any $\a > 0$ and $h\in[H]$,
\begin{align}
\label{eq-agent-fixed-point}
x_h^*&=\Pi_{X_h}[x_h^*-\a \nabla_h F_h(x_h^*,x_{-h}^*)]\\
&=\textstyle\Pi_{X_h}\left[x_h^*- \frac{\a}{N_h}\sum_{i\in \V_h} \nabla_h f^i(x_h^*,x_{-h}^*)\right] .\nonumber
\end{align}
This expression can be rewritten compactly as:
\begin{align}
\label{eq-fixed-point}
x^*=\Pi_{X} [x^*-\a M(x^*) ].
\end{align}

\begin{assumption}\label{assum:map_monotone}
The game mapping $M(\cdot)$ in~\eqref{eq:gamemapping} is strongly monotone with the constant $\mu>0$.
\end{assumption}

\begin{remark}
Assumption~\ref{assum:map_monotone} implies strong convexity of each cluster's cost function $F_h(x_h,x_{-h})$ on $\re^{p_h}$ for every $x_{-h}\in \re^{p-p_h}$ with the constant $\mu$ (cf. Remark 1 of \cite{Tatarenko2021}). The existence and uniqueness of an NE for the game $\Gamma$ is also guaranteed (cf. Theorem 2.3.3 of \cite{FacchineiPang}). The NE can be alternatively expressed as the fixed point solution, as in \eqref{eq-fixed-point}.
\end{remark}

\subsection{Communication Networks}
Consider the partial-decision information scenario, where there is no central coordinator, and agents are restricted to exchanging information solely through peer-to-peer communication.
The communication framework within the multi-cluster game $\Gamma$ is structured into two separate layers: The first layer represents intra-cluster interactions, facilitating communication within the same cluster, without any connection to agents in other clusters. The second layer represents inter-cluster interactions, facilitating global communication among agents irrespective of their cluster affiliation.

\subsubsection{Intra-cluster Interactions}
The interaction among agents within each cluster $h \in [H]$ is represented by a directed graph $\G_h=(\V_h,\E_h)$, specified by the set of edges $\E_h\!\subseteq \V_h\times\V_h$ of ordered pairs of nodes. Associated with $\G_h$ is a weight matrix $\C_h \in \mathbb{R}^{N_h \times N_h}$ that is compliant with the graph $\G_h$, i.e., 
\begin{align*}
\begin{cases}
[\C_h]_{ij}>0, \quad\text{when} (j,i) \in \E_h,\\
[\C_h]_{ij}=0, \quad\text{otherwise.}
\end{cases}
\end{align*} 
Here, each link $(j,i)\!\in\!\E_h$ indicates that agent $i$ receives information from agent $j$ within the same cluster $\V_h$.

\begin{assumption}\label{asum-graph-Gh}
For every cluster $h \in [H]$, the graph $\G_h$ is strongly connected with a self-loop at every node $i \in \V_h$. The weight matrix $\C_h$ is column-stochastic, i.e., $\one^\T \C_h = \one^\T$.
\end{assumption}

\subsubsection{Inter-cluster Interactions}
Interactions between clusters are facilitated by a global communication network represented by the directed graph $\G=(\V,\E)$, where $\E\subseteq \V\times\V$ specifies the set of edges comprising ordered pairs of nodes. This network connects all $N$ agents, enabling inter-cluster communication, with each link $(j,i)$ in $\E$ indicating that agent $i$ receives information from agent $j$ in the game.

The graph $\G$ is associated with a weight matrix $\R \in \re^{N \times N}$, adhering to the connectivity structure of the graph: 
$[\R]_{ij}>0$ when $(j,i) \in \E$ and $[\R]_{ij}=0$ otherwise.

\begin{assumption}\label{asum-graph-G}
The graph $\G$ is strongly connected and has a self-loop at every node $i \in \V$. The weight matrix $\R$ is row-stochastic, i.e., $\R\one = \one$.
\end{assumption}

\begin{remark}
The inter-cluster communication, as represented by the graph $\mathcal{G}$, presents a broader context than scenarios where communication occurs solely between a leader or representative agent from each cluster \cite{MENG2023110919}, as discussed in \cite{Zimmermann2021}. This communication framework allows multiple agents to interact with those outside their respective clusters. If the graph is structured to facilitate communication solely through one agent to other clusters, it reverts to the leader-follower framework outlined in \cite{MENG2023110919}.
\end{remark}

\subsection{Partial-decision Information Notations} \label{sec:Partial-decision-Notations}

In the scenario with partial information, the local cost function $f^i(x_h,x_{-h})$ is exclusive to agent $i \in \V_h$, while the strategy set $X_h$ is only known to agents within that cluster. To navigate privacy constraints and compute the gradient of the local cost function, at each time $k \ge 0$, each agent $i \in \V$ maintains a local variable, as follows,
\begin{align}
    z^{i}(k)= \col( (z_{h}^{i})_{h \in [H]})\in\re^p,
\end{align}
to estimate the strategy of all clusters. Here, $z_{h}^{i}(k)\in\re^{p_h}$ represents the estimate of agent $i$ regarding the decision $x_{h} \in X_h$ of cluster $h$. The estimate of agent $i$ regarding the strategy of all clusters without the $h$-th component is 
\[z_{-h}^{i}(k) = \col( (z_{\ell}^{i})_{\ell \in [H]\backslash \{h\}}) \in \re^{p-p_{h}}.\]

For a solution $x^*$ to be an NE in accordance with Definition~\ref{def-NE}, consensus among agents concerning the local estimate is necessary. Additionally, all estimates should converge to the unique NE $x^*$. Specifically,
\begin{align}\label{eq:goal}
\lim_{k\to \infty} z^{i}(k) = \lim_{k\to \infty} z^{j}(k) = x^*, \text{ for all } i,j \in \V.
\end{align}

The estimates of all $N$ agents in the game can be concatenated to form an estimate matrix at time $k$, denoted as
\begin{align}
    \bz(k)=[z^1(k),\ldots,z^N(k)]^\T \in \re^{N \times p}.
\end{align}
The estimate matrix $\bz_{h}(k)$ contains the estimates of all agents regarding the decision of cluster $h$, while $\bz_{-h}(k)$ excludes the $h$-th component, 
respectively represented as:
\begin{align}
\bz_{h}(k)&=[z_{h}^1(k),\ldots,z_{h}^N(k)]^\T \in \re^{N \times p_h},\\
\bz_{-h}(k)&=[z_{-h}^1(k),\ldots,z_{-h}^N(k)]^\T \in \re^{N \times (p-p_h)}.
\end{align}

Inspired by the gradient-tracking technique outlined in \cite{Angelia2022AB,Nguyen2023AccAB}, each agent $i \in \V_h$ also maintains an auxiliary variable $y_h^{i}(k) \in \re^{p_h}$ to track the average gradients of cluster $h$ at time $k$. The gradient-tracking vector of cluster $h$, encompassing all local gradient-tracking variables of its agents, is defined as:
\begin{align}
    \by_h(k)=[y_h^1(k),\ldots,y_h^{N_h}(k)]^\T \in \re^{N_h \times p_h}.
\end{align}
The gradient-tracking matrix, comprising of these vectors for all clusters arranged in a block diagonal form, is denoted as:
\begin{align}
\bY(k)=\Diag((\by_h(k))_{h\in [H]}) \in \re^{N \times p}.
\end{align}


\section{Algorithm} \label{sec:algo} 
Inspired by the Projected Push-Pull algorithm for distributed optimization \cite{Orhan2023} adapted to the multi-cluster game scenario, we present a distributed algorithm that combines a consensus-based approach with a gradient-tracking technique. This algorithm respects constraints on agents' information access, imposed by the directed intra- and inter-cluster communications. 

Every agent $i\in \V_h$ within cluster $h$ initializes with an arbitrary decision vectors $v^{i}(0)=z^{i}(0)$, where $z_h^{i}(0)\in X_h$ and $z_{-h}^{i}(0)\in \re^{p-p_h}$, and the gradient tracking vector $y_h^i(0) = \nabla_hf^i(v^i(0)) \in \re^{p_h}$. At each time $k$, agents exclusively share scaled gradient tracking information within their own cluster, while the inter-cluster communication network facilitates the exchange of decision estimates. For every time step $k=0,1,\ldots,$ each agent updates its local estimate and the gradient tracking variable according to the following procedure:
\begin{itemize}
\item Perform consensus update: For all $i \in \V$,
\begin{align}\label{alg:consensus}
v^i(k+1) &= \textstyle\sum_{j \in \V} [\R]_{ij} z^{j}(k),
\end{align}
\item Update gradient tracking: For all $h\in [H], i \in \V_h$,
\begin{align}\label{alg:gradient-tracking}
y_h^i(k+1) =& \textstyle \sum_{j \in \V_h} [\C_h]_{ij} y_h^{j}(k)  \\
&+ \nabla_h f^i(v^i(k+1))- \nabla_h f^i(v^i(k)), \qquad \nonumber
\end{align}
\item Update decision estimate: For all $h\in [H], i \in \V_h$,
\begin{subequations}\label{alg:decision}
\begin{align} 
z_h^i(k+1)~ =& ~(1-\gamma) \Pi_{X_h}\left[v_h^i(k+1) \right]  \label{alg:decisionh}\\
&+  \gamma\Pi_{X_h}\left[v_h^i(k+1) - \a y_h^i(k+1) \right],   \qquad \nonumber\\ 
 z_{-h}^i(k+1) =& ~v_{-h}^i(k+1). \label{alg:decision-h}
\end{align}  
\end{subequations}
\end{itemize}
Here, $\a>0$ represents the step-size and $\g\in (0,1)$ is an averaging parameter. In \eqref{alg:consensus}, agents aggregate their own estimates with those of their immediate neighbors through a weighted average, where the weights are determined by the entries of the weight matrix $\R$, associated with the inter-cluster graph $\G$. Subsequently, agents perform a gradient step in the direction of their cluster's respective gradients, using information from the gradient-tracking variable computed with the latest local estimates, as detailed in the update \eqref{alg:decisionh}. Notably, the estimates of decisions from other clusters remain unaltered, as agents lack influence over the decision strategies of other clusters, as described in the update \eqref{alg:decision-h}.

\begin{remark}
Since the action set $X_h$ is known only to agents within cluster $h$, the estimates made by agents regarding the decisions of other clusters might not belong to the respective cluster's action set. Thus, to ensure feasibility, we apply the projection to both terms of the averaging procedure in \eqref{alg:decisionh}. This differs from the algorithm in \cite{Orhan2023} for distributed optimization, which is indeed a special case where all agents belong to the same cluster and are aware of the action set.
\end{remark}

\subsection{Cluster Gradient Tracking}
Equation \eqref{alg:gradient-tracking} governs the update of the cluster gradient tracking variable. Initially, agent $i\in \V_h$ computes a weighted average of their own gradient estimate and those of their neighbors within the same cluster $h$, utilizing the weights specified by the weight matrix $\C_h$ associated with the intra-cluster communication graph $\G_h$. Subsequently, the gradient of the local cost functions is evaluated at the estimates $v^i(k+1)$ and $v^i(k)$, and the difference is integrated into the results of the intra-cluster gradient consensus step. 

Under Assumption~\ref{asum-graph-Gh}, the update in \eqref{alg:gradient-tracking} guarantees the gradient tracking property at each time $k \ge 0$. Specifically,

\begin{lemma}[\cite{Angelia2022AB}, Lemma 4] \label{lem-sumgrad}
Let Assumption~\ref{asum-graph-Gh} hold. For all $h \in [H]$, 
we have
\[\textstyle\sum_{i\in \V_h} y_h^i(k)=\sum_{i\in \V_h} \nabla_h f^i(v^i(k)), \text{for all} k\ge 0.\]
\end{lemma}


\subsection{Compact Form}
To express the algorithm in compact form, we introduce the assignment matrix, defined as:
\begin{align} \label{eq-Qmatrix}
\bQ_h^i &= [\zero_{p_h \times p_{<h}},\mbi_{p_h},\zero_{p_h \times p_{>h}}] \in \re^{p_h \times p},
\end{align}
where $p_{<h} = \sum_{\ell=1}^{h-1} p_{\ell}$ and $p_{>h} = \sum_{\ell=h+1}^{H} p_{\ell}$. The matrix $\bQ_h^i$ selects agent $i$'s estimate of their own cluster's decision in a stacked vector $z^i$, thus, $\bQ_h^i z^i = z_h^i \in \re^{p_h}$.

Define the set 
\[\Omega = \{\bz \in \re^{N\times p}| \bQ_h^i \bz_{[i,:]}^\T \in X_h, \forall i\in \V_h, h\in [H]\},\]
where $\bz_{[i,:]}$ represents row $i$ of the matrix $\bz$. This set guarantees that the estimate of agent $i\in \V_h$ regarding cluster $h$'s decision is within the strategy set $X_h$, while their estimate concerning the decision of other clusters is in $\re^{p-p_h}$.

We define the matrix $\bG_h(k)$ containing the gradients of cluster $h$ evaluated using the aggregated estimates $\bv(k)$ as:
\[\bG_h(k) = \col\left( (\nabla_h f^i(\bv(k))^\T)_{i \in \V_h}\right).\]

Denoting related notations for the variable $\bv \in \re^{N\times p}$ using similar notation conventions as those introduced in Section~\ref{sec:Partial-decision-Notations} for the variable $\bz$, we can represent the updates in \eqref{alg:consensus}--\eqref{alg:decision} in compact forms as follows:
\begin{subequations}\label{alg-comp}
\begin{align} 
\!\!\!\bv(k+1) &\!=\! \R \bz(k)\label{alg-comp-avg}\\
\!\!\!\by_h \!(k\!+\!1) &\!=\! \C_h\by_h(k) \!+\! \bG_h(\bv(k+1)) \!-\! \bG_h(\bv(k)), \forall h \!\!\label{alg-comp-estimate}\\
\!\!\!\bz(k+1) &\!=\! (1-\g)\Pi_{\Omega} [\bv(k+1)] \nonumber\\
&\quad + \g\Pi_{\Omega} [\bv(k+1) - \a\bY(k+1)]. \label{alg-comp-grad}
\end{align}  
\end{subequations}

\section{Analysis}\label{sec:conv_results}
\subsection{Preliminaries}\label{sec:basicre}
\noindent We have the following results regarding the weight matrices:

\begin{lemma}[\cite{pshi21}, Lemma 1] \label{lem-ChRhEig}
Let  $\{\C_h\}_{h \in [H]}$ satisfy Assumption~\ref{asum-graph-Gh}, and let $\R$ satisfy Assumption~\ref{asum-graph-G}. It follows that:
\begin{itemize}
\item[(i)] Each matrix $\C_h$ has a unique, positive right eigenvector $\pi_h$ corresponding to eigenvalue $1$, i.e., $\C_h\pi_h = \pi_h$, such that $\one^\T \pi_h = 1$, for all $h \in [H]$.
\item[(ii)] The matrix $\R$ has a unique, positive left eigenvector $\phi$ corresponding to eigenvalue $1$, i.e., $\phi^\T \R = \phi^\T$, such that $\phi^\T \one = 1$.
\end{itemize}
\end{lemma}

\begin{lemma}[\cite{Angelia2022AB}, Lemma 5.5]
\label{lem-Contraction-Col}
Consider $\{\C_h\}_{h \in [H]}$ and let Assumption~\ref{asum-graph-Gh} hold. For arbitrary $\bu \in \re^{N_h \times p_h}$, we have:
\begin{align}
    \|\C_h\bu - \pi_h\one^\T\bu\|_{\pi_h^{-1}} &\leq \sigma_C\|\bu - \pi_h\one^\T\bu\|_{\pi_h^{-1}}.
\end{align}
Here, $\sigma_C = \max_{h\in[H]}\left\{\!\sqrt{1\!-\!\frac{\min^2(\pi_h)\left({\min}^{+}(\C_h)\right)^2}{\max^3(\pi_h)\mathsf{D}(\G_h)\mathsf{K}(\G_h)}}\right\} \!\in\! (0,1).$
\end{lemma}

\begin{lemma}[\cite{nguyen2022distributed}, Lemma 6] \label{lem-Contraction-Row}
Consider the matrix $\R$ and let Assumption~\ref{asum-graph-G} hold. For arbitrary $\bu, \bv \in \re^{N \times p}$, we have:
\begin{align}\label{eq-lem-Contraction-Row-eq1}
\|\R\bu - \bv\|_{\phi}^2 \leq \|\bu - \bv\|_{\phi}^2 - c\|\bu - \one\phi^\T\bu\|_{\phi}^2.
\end{align}
where $c^2=\frac{\min(\phi)\left({\min}^{+}(\R)\right)^2}{\max^2(\phi)\mathsf{D}(\G)\mathsf{K}(\G)}$. Additionally, choosing $\bv = \one\phi^\T\bu$ and defining $\sigma_R = \sqrt{1-c^2} \in (0,1)$ results in
\begin{align}
\|\R\bu - \one\phi^\T\bu\|_{\phi} &\leq \sigma_R\|\bu - \one\phi^\T\bu\|_{\phi}.
\end{align}
\end{lemma}

\begin{lemma}[\!\! \cite{Tatarenko2021}, Lemma 1] \label{lemma:LipschitzNablafi}
Let  Assumption~\ref{assum:lip}(iii)--(iv) hold. 
For all $x,y \in \re^p$, 
we have for all $h \in [H]$ and $i \in \V_h$:
\begin{align}\label{eq-DelJLipschitz}
\|\nabla_h f^i(x)-\nabla_h f^i(y)\|^2\le L^2\|x-y\|^2,
\end{align}
where $L = \max_{i\in\V} \sqrt{(L_1^i)^2+(L_2^i)^2}$.
\end{lemma}


\subsection{Main Results}\label{sec:conv}


The convergence of the proposed algorithm is analyzed based on establishing a contraction relationship among three critical metrics:
(i) the optimality gap $\|\bz(k) - \bx^*\|_{\phi}$, where $\bx^* = \one_N({x^*})^\top \in \re^{N\times p}$; 
(ii) the consensus error $\bD(\bv(k))$
\begin{align*}
\bD(\bv(k)) = \sqrt{\sum_{i \in \V_h}\sum_{j \in \V_h} \phi_i\phi_j\|v^i(k) - v^j(k)\|^2}; 
\end{align*}
and, (iii) the gradient tracking error \begin{align*}
\bS(\by(k)) &= \textstyle\sum_{h\in [H]} \bS_h(\by(k)),
\end{align*}
where for all $h \in [H]$:
\begin{align*}
\bS_h(\by(k)) &=\! \sqrt{\sum_{i\in \V_h}\!\frac{1}{[\pi_h]_i}\bigg\|y_h^i(k) - [\pi_h]_i \!\sum_{j\in\V_h} y_h^j(k)\bigg\|^2}.
\end{align*}

The contraction relation for $\|\bz(k) - \bx^*\|_{\phi}$ is presented 
below, and it is the cornerstone of our analysis.
\begin{lemma}\label{re-OptimalityGap}
Let Assumption~\ref{assum:lip}--\ref{asum-graph-G} hold. There exists a positive constant $\bar{\a}_Q$ such that for  $\a \in (0,\bar{\a}_Q)$ and  $\g \in (0,1)$, there exists $\rho_{\a}\in(0,1)$ such that:
\begin{align*}
\|\bz(k)-\bx^*\|_{\phi} &\le (1-\g \rho_{\a})\|\bz(k-1)-\bx^*\|_{\phi} \\
 &+ \a\g\varphi\sqrt{N}L \bD(\bv(k)) + \a\g \bS(\by(k)),
\end{align*}
for all $k \ge 1$, where $\varphi = \sqrt{1/\min(\phi)}$.
\end{lemma}


\begin{proof}
To characterize the scenario where agents have information about cluster gradients, we define
for all $h\in[H]$ and $i\in\V_h$: 
\begin{align*}
&w_h^i(k) = (1-\g) \Pi_{X_h}[v_h^i(k)]  \\
&\qquad\qquad +\g\Pi_{X_h}[v_h^i(k)- \a N_h[\pi_h]_i \nabla_h F_h(v^i(k))],\nonumber\\
&w_{-h}^i(k) = v_{-h}^i(k).
\end{align*}

We define the scaled gradient matrix $\tilde\bG(\bv(k))$ as:
\begin{align*}
\tilde\bG(\bv(k)) = \Diag((\tilde\bG_h(\bv(k))_{h\in [H]}) \in \re^{N \times p},
\end{align*}
where $\tilde\bG_h(\bv(k)) = \col\left( (N_h[\pi_h]_i\nabla_h F_h(v^i(k))^\T)_{i \in \V_h}\right)$.
Then, the compact form for $\bw(k)$ can be obtained as follows:
\begin{align}
\bw(k) = (1-\g)\Pi_{\Omega} [\bv(k)] + \g\Pi_{\Omega} [\bv(k) - \a\tilde\bG(\bv(k))]. \label{eq-wk-comp}
\end{align}  
By the triangle inequality, we obtain
\begin{align}\label{eq:addw}
\|\bz(k)-\bx^*\|_{\phi} \le\|\bw(k)-\bx^*\|_{\phi} + \|\bz(k)-\bw(k)\|_{\phi}.
\end{align}
For the first term in \eqref{eq:addw}, using \eqref{eq-wk-comp}, the fact that $x^* \in X$, the non-expansiveness property of the projection yields
\begin{align}\label{eq:wxstar}
\|\bw(k)-\bx^*\|_{\phi} &\le (1-\g)\|\bv(k)-\bx^*\|_{\phi}\\
&+\g\|\bv(k) - \bx^* -\a(\tilde\bG(\bv(k))-\tilde\bG(\bx^*))\|_{\phi}.\nonumber
\end{align}

Following the proof of Lemma~3 in \cite{bianchi2024end} with $\bv(k) = \R \bz(k-1)$, there exists a positive upper bound $\bar{\a}_Q$ such that for any step-size $\a \in (0,\bar{\a}_Q)$, we can derive the matrix $\bar{Q}_\a$ with the largest eigenvalue $\bar{\lambda}(\bar{Q}_\a) \in (0,1)$, and the next relation holds:
\begin{align}\label{eq:wxstar1}
&\|\bv(k) - \bx^* -\a(\tilde\bG(\bv(k))-\tilde\bG(\bx^*))\|_{\phi} \\
\le& \sqrt{\bar{\lambda}(\bar{Q}_\a)} \|\bz(k\!-\!1) - \bx^*\|_{\phi} = (1-\rho_{\a})\|\bz(k\!-\!1) - \bx^*\|_{\phi},\nonumber
\end{align}
where $\rho_{\a}=1-\sqrt{\bar{\lambda}(\bar{Q}_\a)} \in (0,1)$.

Using the update in \eqref{alg-comp-avg} and by applying relation \eqref{eq-lem-Contraction-Row-eq1} in Lemma~\ref{lem-Contraction-Row}, with $\bu = \bz(k-1)$ and $\bv = \bx^*$, we obtain:
\begin{align}\label{eq-vzrel}
\!\!\|\bv(k)\!-\!\bx^*\|_{\phi} \!=\! \|\R\bz(k\!-\!1)\!-\!\bx^*\|_{\phi} \!\le\! \|\bz(k\!-\!1)\!-\!\bx^*\|_{\phi}.\!\!\!
\end{align}
Combining the preceding relation with \eqref{eq:wxstar} and \eqref{eq:wxstar1} yields
\begin{align}\label{eq:wxstar4}
\|\bw(k)-\bx^*\|_{\phi} \le (1-\g \rho_{\a})\|\bz(k-1)-\bx^*\|_{\phi}.
\end{align}

For the second term in \eqref{eq:addw}, we have
\begin{align*}
\!\!\|\bz(k)-\bw(k)\|_{\phi}^2 
= \frac{1}{H}\!\!\sum_{h\in [H]}\sum_{i\in \V_h}\![\phi_h]_i \|z_h^i(k)-w_h^i(k)\|^2.\!\!
\end{align*}
Using the update for $z_h^i(k)$ in \eqref{alg:decisionh} and the non-expansiveness property of the projection, we obtain
\begin{align*}
\|z_h^i(k)-w_h^i(k)\| \le & \a\g\|y_h^i(k) - N_h[\pi_h]_i \nabla_h F_h(v^i(k))\|\\
\le & \a\g \bigg\|y_h^i(k) - [\pi_h]_i \sum_{j\in \V_h} \nabla_h f^j(v^i(k))\bigg\|.
\end{align*}
Using Lemma~\ref{lem-sumgrad} and the triangle inequality, we obtain
\begin{align}\label{eq:zw2}
&\|\bz(k)-\bw(k)\|_{\phi} \nonumber\\
\le & \a\g \sqrt{\sum_{h\in [H]}\sum_{i\in \V_h}\phi_i\bigg\|y_h^i(k) - [\pi_h]_i\sum_{j\in\V_h} y_h^j(k)\bigg\|^2} \nonumber\\
+& \a\g\sqrt{N}\!\sqrt{\sum_{\substack{h\in [H]\\i,j\in \V_h}}\!\phi_i[\pi_h]_i^2 \bigg\| \nabla_h f^j(v^j(k)) -  \nabla_h f^j(v^i(k))\bigg\|^2}\nonumber\\
\le & \a\g \!\!\!\sum_{h\in [H]} \!\!S_h(\by(k)) + \a\g\sqrt{N} L \!\!\sqrt{\!\sum_{\substack{h\in [H]\\i,j\in \V_h}} \!\phi_i \| v^j(k) -  v^i(k)\|^2}\nonumber\\
\le & \a\g \bS(\by(k)) +\a\g\varphi\sqrt{N}L \bD(\bv(k)).
\end{align}
Here, we use Lemma~\ref{lemma:LipschitzNablafi}, and the fact that $\phi$ and $\pi_h$ are stochastic vectors.
Then, the result follows directly from \eqref{eq:addw}, \eqref{eq:wxstar4} and \eqref{eq:zw2}.
\end{proof}

The contraction relation for $\bS(\by(k))$ follows analogous reasoning to the analysis presented in Proposition 5 of \cite{Orhan2023} and Proposition 1 of \cite{Zimmermann2021}. Similarly, the contraction relation for $\bD(\bv(k))$ follows along the lines of Lemma 7 of \cite{Orhan2023}. These analyses leverage the inherent row- and column-stochasticity properties of the matrices $\R$ and $\C_h$ in Assumption~\ref{asum-graph-G} and Assumption~\ref{asum-graph-Gh}, respectively, in conjunction with Lemma~\ref{lem-Contraction-Col}, Lemma~\ref{lem-Contraction-Row} and relation \eqref{eq-vzrel}. The contraction coefficients take the form 
$\sigma + \a\g\psi$, where $\sigma<1$ represents the contraction coefficient in Lemma~\ref{lem-Contraction-Col} and Lemma~\ref{lem-Contraction-Row}, and $\psi$ is a constant related to the properties of the game $\Gamma$ and the communication network; importantly, $\lim_{(\a,\g)\to\zero} \a\g = 0$. Consequently, there exists a choice of the step-sizes $\a>0$ and $\g>0$ such that $\sigma + \a\g\psi < 1$, which is crucial in controlling the spectral radius of the coefficient matrix in the composite relation of the errors, ensuring convergence. 

In particular, by defining the error vector $\bE(k)$ as
\begin{align*}
\bE(k)=\Big(\|\bz(k) - \bx^*\|_{\phi},\bD(\bv(k)),\bS(\by(k))\Big)^{\T},
\end{align*}
we have the following result:
\begin{theorem}
Let Assumption~\ref{assum:lip}--\ref{asum-graph-G} hold. There exists a positive constant $\bar{\a}_{\max}\le \bar{\a}_Q$ such that for  $\a \in (0,\bar{\a}_{\max})$, there exists $0<\bar{\g}_{\max}^{\a}<1$ dependent on $\a$, such that for $\g \in (0,\bar{\g}_{\max}^{\a})$, a composite relation for the errors can be established as follows:
\[\bE(k+1)\le \textrm{M}(\a,\g)\bE(k),\qquad \hbox{for all ~~ $k\ge0$},\]
where the matrix $\textrm{M}(\alpha,\gamma)$ 
 takes the form
\begin{align*}
\!\!\!\left[\!\!
\begin{array}{ccc}
	 1 -\g \rho_{\a} 
   + \mc{O}(\a\g^2)
  &  \mc{O}(\a\g)  & \mc{O}(\a\g)  \cr
	 \mc{O}(\g) & \sigma_R + \mc{O}(\a\g)   & \mc{O}(\a\g) \cr
	 \mc{O}(\g) & \mc{O}(1) & \sigma_C +\mc{O}(\a\g) 
	 \end{array}\!\!\right]\!\!\!\!\!
\end{align*}
and has a spectral radius smaller than $1$. 
Therefore, $\lim_{k\to\infty} \|z_h^i(k) -x^*\|=0$ with a linear convergence rate of the order $\mathcal{O}\Big(\rho_M^k\Big)$, for all $h\in[H]$ and $i\in \V_h$,
with $\rho_{M}<1$ being the spectral radius of $M(\a,\g)$.
\end{theorem}


Given that $\rho_{\a}\in (0,1)$, $\sigma_R\in (0,1)$ and $\sigma_C\in (0,1)$, selecting suitable step-sizes $\a \in (0,\bar{\a}_{\max})$ and $\g \in (0,\bar{\g}_{\max}^{\a})$ ensures that all diagonal entries of $\textrm{M}(\a,\g)$ are less than $1$ and $\det(\mbi-\textrm{M}(\a,\g))>0$, thereby guaranteeing $\rho_{\textrm{M}}<1$ (cf. Lemma~8 in \cite{pshi21}).
Consequently, the proposed algorithm exhibits a linear convergence rate. This procedure mirrors \cite[Propositions~1 and Theorem~1]{Orhan2023},
among others (cf. \cite{Angelia2022AB,Nguyen2023AccAB,Zimmermann2021,Zimmermann2022,pang2023distributed}), which we omit here.

\begin{remark}
The convergence of algorithm \eqref{alg:consensus}--\eqref{alg:decision} is ensured for sufficiently small step-sizes $\a$ and $\gamma$, and explicit bounds could also be computed. Yet,  practical applications often requires manual optimization of step-sizes, given the conservative nature of theoretical bounds as also noted in prior research (cf. \cite{nguyen2022distributed,pshi21}).
\end{remark}

\begin{remark} [Impact of Averaging Parameter $\g$] ~\sloppy 
The averaging procedure in \eqref{alg-comp-grad}, governed by the parameter $\g$, is crucial for the algorithm's convergence. We refer to the proof provided for distributed optimization in Section VI.D. of \cite{Orhan2023} to demonstrate that even for the special case of the game $\Gamma$ when there is only one cluster ($H=1$), the averaging procedure remains essential. For completeness, we restate the statement herein. 

Suppose we eliminate the averaging step by setting $\g = 1$ in equation~\ref{alg-comp-grad}, resulting in the update rule:
\[\bz(k+1)  = \Pi_{\Omega} [\bv(k+1) - \a\bY(k+1)].\]
Section VI.D. of \cite{Orhan2023} establishes that with this update rule, it becomes infeasible to constrain the following expression:
\begin{align*}
&\|\bv(k+1)-\bv(k)\| \\
\le & c_1(\a)\|\bz(k)-\bx^*\|_{\phi} + c_2(\a)\bD(\bv(k)) + c_3(\a) \bS(\by(k)),
\end{align*}
such that $\lim_{\a\to 0} c_1(\a) = 0$. The term $\|\bv(k+1)-\bv(k)\|$ appears when establishing the contraction relation for the gradient tracking error $\bS(\bv(k+1))$. This error is pivotal for establishing the composite relations, which, in turn, are crucial for proving the convergence of the algorithm, as detailed in Section~\ref{sec:conv}.

\end{remark}


\section{Energy Management in Networked Microgrids} \label{sec:simulation}
\subsection{System Model}
\begin{figure}[hbt!]
\centering
\includegraphics[width=0.45\textwidth,height=0.16\textheight]{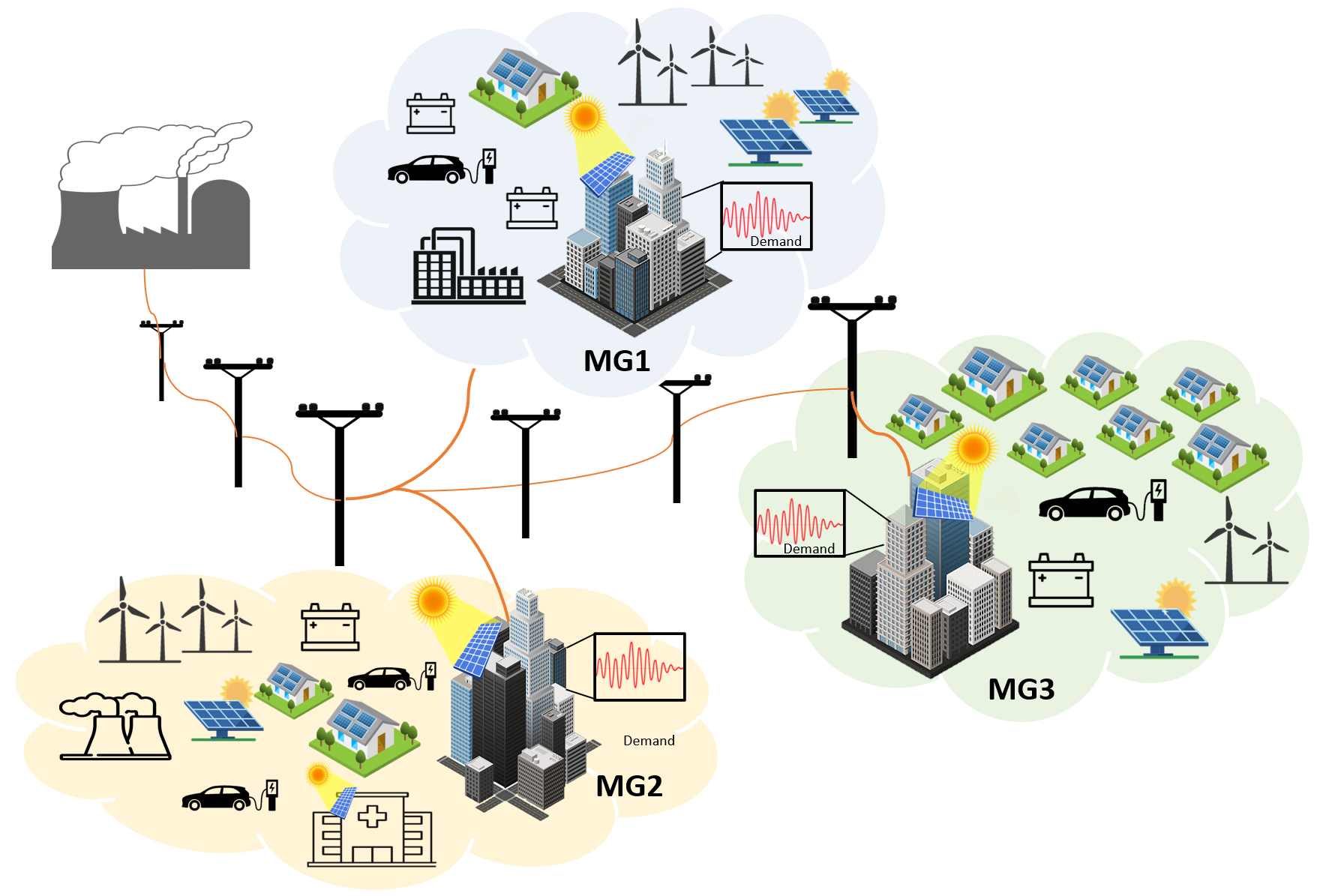}
\caption{System model}
\label{fig:model}
\vspace{-0.2cm}
\end{figure}

In a day-ahead energy management problem, in line with references \cite{Liu2021,Belgioioso2020,Zimmermann2022}, we consider a system of $H$ microgrids (MGs) or distinct energy systems, each tasked with supplying power to its consumers over a predefined time horizon $\mc{T}=\{1, \ldots, T\}$. The system consists of a set $\V$ of $N$ components. MG $h$ is equipped with $N_h^g$ energy generation units and $N_h^b$ battery units. In total, the MG system comprises $N$ components, namely $\sum_{h\in [H]}(N_h^g+N_h^b)=N$. The sets $\V_h^g$ and $\V_h^b$ represent the respective partitions of the set $\V$ pertaining to energy generation and storage components within MG $h$. At each time slot $t \in \mc{T}$, MG $h \in [H]$ aims to meet the power demand $\PD_h(t)$ of its consumers while minimizing the associated cost of power provision, which includes operational expenses and grid procurement. The system model is depicted in Fig~\ref{fig:model}.

\subsection{Problem Formulation}
Let $\PR^i(t)$ represent the power generated by generator $i$ at time $t$, subject to the capacity limits:
\begin{align}\label{eq-num-GenConsts}
\PR^{i,\min} \le \PR^i(t) \le \PR^{i,\max}, \text{for all} i \in \V_h^g.
\end{align}
The cost incurred by generator $i \in \V_h^g$ is given by \cite{Liu2021}:
\begin{align*}
\FR^i(\PR^i(t)) = a^i (\PR^i(t))^2 + b^i \PR^i(t) + c^i,
\end{align*}
where $a^i$, $b^i$, and $c^i$ are positive constants, for all $i \in \V_h^g$. 
Let $\PR_h(t) = \col\big((\PR^i(t))_{i\in\V^g_h}\big)$, we have the constraint set:
\begin{align*}
    \Omega_h^{\PR} = \{\PR_h(t)| \PR^i(t) \in \re^+, \eqref{eq-num-GenConsts}, \forall i\in\V^g_h, \forall t \in \mc{T}\}.
\end{align*}

Let $\PB^i(t)$ denote the power flow in battery $i$ at time $t$. A positive value of $\PB^i(t)$ indicates discharge by an amount of $|\PB^i(t)|$, while a negative value signifies charging. The utilization of battery $i\in \V_h^b$ results in a penalty \cite{Belgioioso2020}:
\begin{align*}
\FB^i(\PB^i(t)) = a^i (\PB^i(t))^2 + b^i |\PB^i(t)| + c^i,\!\!
\end{align*}
where $a^i$, $b^i$, and $c^i$ are positive constants, for all $i \in \V_h^b$. The amount of power flow in the battery is constrained by:
\begin{align}\label{eq-num-BatteryConsts}
    \PB^{i,\min} \leq \PB^i(t) \leq \PB^{i,\max}.
\end{align}

The battery's charging or discharging capability is determined by its charge $\PC^i(t)\in[0,\PC^{i,\max}]$, we have \cite{Belgioioso2020}:
\begin{align*}
    -\PC^i(t) \leq - \PB^i(t) \leq \PC^{i,\max} - \PC^i(t).
\end{align*}
Let $\eta^i\in(0,1)$ denote energy loss over time when the battery is idle, yielding $\PC^i(t+1) = \eta^i \PC^i(t)$. The charge at the start of time slot $t$ can be computed as 
\begin{align*}
\textstyle  \PC^i(t) = (\eta^i)^{t-1} \PC^i(1) - \sum_{s=1}^{t-1} (\eta^i)^{t-s} \PB^i(s).
\end{align*}
Consequently, we have the following constraint:
\begin{subequations}\label{eq-num-BatteryConst2}
 \begin{align}\label{eq-num-BatteryConst2a}
-(\eta^i)^{t-1} \PC^i(1) &\leq - \textstyle \sum_{s=1}^{t} (\eta^i)^{t-s} \PB^i(s),\\
\textstyle - \sum_{s=1}^{t} (\eta^i)^{t-s} \PB^i(s) &  \leq \PC^{i,\max} - (\eta^i)^{t-1} \PC^i(1). \label{eq-num-BatteryConst2b}
\end{align}   
\end{subequations}

Furthermore, by the end of the planning horizon, the battery charge $\PC^i(T)$ approaches the desired state of charge $\PC^{i,\textnormal{des}}$ \cite{Atzeni2013}, i.e., for some small $\epsilon_i > 0$:
\begin{align}\label{eq-num-BatteryConste}
    |\PC^i(T) - \PC^{i,\textnormal{des}}| \leq \epsilon_i.
\end{align}

The constraint set for $\PB_h(t) = \col\big((\PB^i(t))_{i\in\V^g_h}\big)$ is:
\begin{align*}
   \!\!\! \Omega_h^{\PB} = \{\PB_h(t)|\PB^i(t)\in\re, \eqref{eq-num-BatteryConsts}-\eqref{eq-num-BatteryConste}, ~\forall t \in \mc{T}, \forall i \in \V_h^b\}.\!\!\!
\end{align*}

When the generated energy falls short of meeting the power demand $\PD_h(t)$, MG $h$ must acquire power of amount $\PG_h(t)\in [0,\PG_h^{\max}]$ from the market,
at the price \cite{Belgioioso2020},
\begin{align*}
p_B(\PG(t)) = \textstyle \zeta\left(\sum_{h=1}^H \PG_h(t)\right),
\end{align*}
where $\PG(t) = \col( (\PG_h(t))_{h \in [H]})$ 
and $\zeta$ is a positive scaling factor. 
We assume that MGs are allowed to sell back surplus energy to the grid when the onsite generation exceeds the demand at a sell-back price of $p_S(\PG(t)) = \varrho p_B(\PG(t))$. The sell-back price is expected to be less than the procurement cost, thus, $\varrho \in (0,1)$.
Let  $\PS_h(t)$ represent the amount of electricity to be sold to the grid from MG $h$ at time $t$. Then, we can calculate the electricity cost as follows:
\begin{align*}
\FE_h(\PG(t),\PS_h(t)) = p_B(\PG(t)) (\PG_h(t) - \varrho \PS_h(t)).
\end{align*}

The optimization problem for MG $h$ is 
\begin{align*}
&~~ \min_{x_h} ~~ \textstyle\sum_{t=1}^T \left[\sum_{i \in \V_h^g}\FR^i(\PR^i(t)) + \sum_{i \in \V_h^b}\FB^i(\PB^i(t))\right]\\
& \qquad \qquad \textstyle+\sum_{t=1}^T\FE_h(\PG(t),\PS_h(t))\\
&\text{s.t.} \PR_h(t) \in \Omega^{\PR}_h,~~ \PB_h(t) \in \Omega^{\PB}_h, ~~\forall t, \\
&\qquad\quad~\! \PG_h(t) \in [0,\PG_h^{\max}],~~  \PS_h(t) \ge 0, ~~\forall t, \\
&\sum_{i \in \V_h^g}\!\PR^i(t) +\!\! \sum_{i \in \V_h^b}\!\PB^i(t)+ \PG_h(t) = \PD_h(t) + \PS_h(t), ~\forall t,\!
\end{align*}
where $x_h = \col\left(\PR_h(t),\PB_h(t),\PG_h(t),\PS_h(t)\right)$.

\subsection{Simulation Results}
We consider $H\!=\!6$ MGs with a total of $N\!=\!50$ units selected randomly as generators or batteries over a time horizon of $T\!=\!24$ hours. Generators are chosen randomly from a mix of traditional (coal, natural gas, oil) and renewable (wind, solar, nuclear, hydropower) sources, with power limits and cost coefficients obtained from the MATPOWER dataset\footnote{\textnormal{https://github.com/MATPOWER/matpower/tree/master/data}}. 
For battery units, cost coefficients are generated as $a^i \!\sim\! U[0.1,5]\$$/MWh$^2$, $b^i \!\sim\! U[5,50]\$$/MWh, and $c^i \!\sim\! U[-50,50] \$$. Battery systems have capacities ranging from $50$MWh to $200$MWh. The battery leakage rate ranges from $0.95$ to $0.99$. The maximum charge rate is randomly selected from the range of $0.8$C to $1$C, where C represents the battery's capacity, while the initial charge ranges from $0.2$C to $0.5$C. The rate of the sell-back price is $\varrho=0.8$. The demand $\PD_h^i(t)$ is randomly generated within the range of $[500,2000]$MWh. The code repository is available\footnote{\textnormal{https://github.com/duongnguyen1601/Distributed\_NCluster\_Game/}}.


The NE $x^*$ is computed using the update in~\eqref{eq-agent-fixed-point}, assuming agents have full information access. In the partial-information decision scenario, we can verify that all the assumptions are satisfied for this problem. Specifically, the cost function takes the form of a quadratic function with respect to the decision variable, and it can be confirmed that the game mapping exhibits strong monotonicity. Additionally, the constraint set of each cluster includes a linear equality constraint, along with the convex and compact sets $\Omega_h^{\PR}$ and $\Omega_h^{\PB}$, thus, it is closed and convex. To ensure strong connectivity among graphs, we establish a directed cycle linking all agents within each cluster and linking all clusters. Thus, our proposed distributed algorithm can effectively be applied to estimate the NE. The optimality gap between the estimates obtained using our approach and the NE $x^*$ is depicted in Figure~\ref{fig:Residual}(a), while the total cost for each MG is depicted in Figure~\ref{fig:Residual}(b), demonstrating the convergence property of the proposed algorithm.

\begin{figure}[ht!]
\centering
\subfigure[Optimality gap (log-scale)]{
	\includegraphics[width=0.23\textwidth]{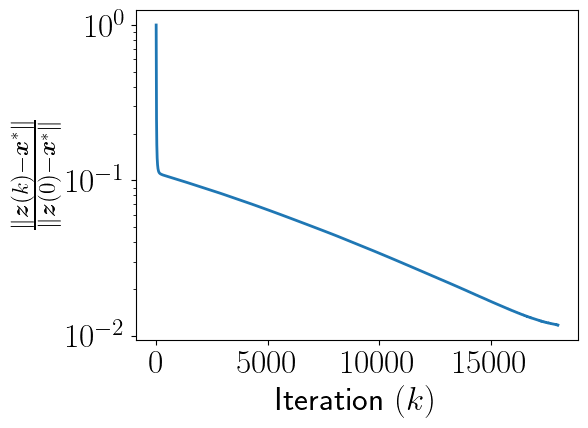}\label{fig:Residual}
	}  
    \hspace*{-.8em}
\subfigure[Total cost]{
	\includegraphics[width=0.21\textwidth]{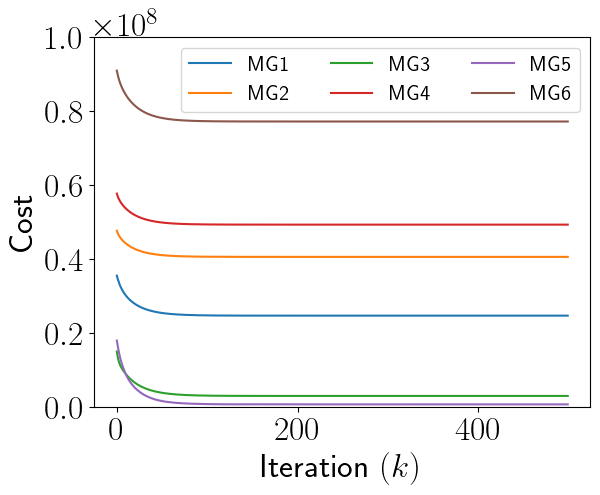}\label{fig:Error}  }
\caption{Convergence plots.}
\label{fig:Residual}
\vspace{-0.3cm}
\end{figure}

\section{Conclusions and Future Work}
\label{sec:conc}
We propose a distributed NE seeking algorithm for multi-cluster games, effectively addressing technical challenges associated with constraint sets and directed communication networks. Our algorithm integrates consensus-based methods, gradient-tracking techniques, and an averaging procedure; the latter is crucial to control the optimality gap and to ensure linear convergence.
We validate the linear convergence of our algorithm through an energy management problem in networked microgrids. In future works, we aim to extend our algorithm to handle more complex scenarios, such as dynamic games and dynamic communication networks.

\bibliographystyle{IEEEtran}
\bibliography{references}

\end{document}